\documentclass[a4,11pt,reqno]{amsart}
\usepackage{amsfonts,amssymb}
\theoremstyle{plain}
\newtheorem{theorem}{Theorem}
\newtheorem{lemma}{Lemma}
\newtheorem{corollary}{Corollary}

\newcommand{\eps}{\varepsilon}

\newcommand{\Om}{\Omega}
\usepackage{amsfonts}
\newcommand{\bproof}{{\bf Proof.}{\hspace{0.2cm}}}
\newcommand{\eproof}{{$\Box$}{\medskip}\\}

\newcommand{\mes}{\rm{mes}}

\newcommand{\supp}{\rm{supp}}
\newcommand{\varep}{\varepsilon}

\newcommand{\R}{\mathbb{R}}
\newcommand{\N}{\mathbb{N}}

\def\cal{\mathcal}

\def\Ge{\mathcal{G}}
\def\E{\mathcal{E}}

\theoremstyle{definition}
\newtheorem{definition}{Definition}

\theoremstyle{remark}
\newtheorem{remark}{Remark}
\newtheorem{example}{Example}
\newcommand{\ssubset}{\subset\subset}
\newcommand{\wf}{{\rm WF}}

\begin{document}
\title{Real analytic generalized functions}
\author{S.\ Pilipovi\'{c}, D.\ Scarpalezos,  V.\ Valmorin}
\address{
Department of Mathematics and Informatics, University of Novi Sad,
Trg D. Obradovi\'ca 4, 21000 Novi Sad, Serbia}
\address{U.F.R.\ de Math\'{e}matiques,
Universit\'e Paris 7,
2 place Jussieu,
Paris 5 $\buildrel e \over =$,
75005, France}
\address{
Universit{\'e} des Antilles et de la
Guyane, D{\'e}partement Math-Info, Campus de Fouillole:
97159 Pointe {\`a} Pitre Cedex, France}
\date{}
\maketitle
\begin{abstract}
Real analytic generalized functions are investigated as well as
the analytic singular support and analytic wave front of a
generalized function in $\Ge(\Om)$ are introduced and described.


\end{abstract}

\section{Introduction and definitions}

Holomorphic generalized functions are defined as the solutions to
$\bar{\partial}_{z_i}G=0, i=1,...,p,$ where $G$ is an element of a
generalized function algebra $\mathcal{G}(\Omega), \Omega
\subset\mathbb{C}^p$ (see \cite{bia}, \cite{col}, \cite{co1},
\cite{gkos},  \cite{ober 001} and \cite{dim} for the definition
and the properties for generalized function algebra).  Recall,
elements of this algebra are classes of equivalences of nets of
smooth functions with respect to an ideal $\mathcal{N}(\Omega)$
which consists of nets $(n_\varepsilon)_\varepsilon$ converging to
zero, on compact sets, faster than any positive power of
$\varepsilon$ as $\varepsilon\rightarrow 0.$ So,
$\bar{\partial}G=0$ means
$\bar{\partial}G_\varepsilon=O(\varepsilon^a)$ on compact sets  of
$\Omega$,
 for any $a>0$. Thus,
various  properties for such generalized functions have became interesting and important.
The algebra $\mathcal{G}_H(\Omega), \Omega \subset\mathbb{C}^p,$ of  holomorphic generalized functions
is  introduced and studied
 in \cite{col}, \cite{ar1}, \cite{coga2} and \cite{coga} while in  \cite{pili} is
 proved the existence of a global holomorphic representative (any member of a
 representative is holomorphic)
 for an $f\in\mathcal{G}_H(\Omega), \Omega \subset \mathbb{C}.$ Moreover, in \cite{kelsca} is shown
(in the one dimensional case) that a holomorphic generalized
 function equals zero if it is equal to zero at every classical point.

 In this paper we are
interested in the class of real analytic generalized functions
$\mathcal{G}_A(\omega),\; \omega$ is open in $\mathbb{R}^d,$ in
order to study generalized  analytic microlocal regularity of
generalized functions.

We give in Section 2 a  representation of an arbitrary generalized
function $f\in \mathcal{G}(\omega)$, $\omega $ is open in
$\mathbb{R}^d,$
as a boundary value on $I\subset\subset \omega$ of a corresponding
holomorphic generalized function $[(F_\eps)_\eps]$ on $I\times
(-1/\eta,1/\eta)^d$ which means
$(\bar{\partial}_{z_i}F_\eps)_\eps\in\mathcal{N}(I\times
(-1/\eta,1/\eta)^d), i=1,...,d$.
 Further, in the case  $d=1$
 we give a global real analytic representative of a real analytic generalized function.
 Then we define
 analytic singular support of a generalized function. Also in Section 2 we prove
that a  generalized real analytic function equals zero if it is
equal to zero at every classical point of a subset of $\omega$ of positive measure by proving
 the same result for a holomorphic generalized  function in an open
set of $\mathbb{C}^p$ (see \cite{kelsca} for the case $p=1$).
Generalized analytic wave front of $f\in\mathcal{G}(\omega)$ is
introduced in Section 3 (see \cite{horkun}, \cite{gar1}, \cite{gar2},
\cite{garhor} and \cite{horman} for the local and microlocal
$\mathcal{G}^\infty$-properties of Colombeau-type generalized
functions). It is proved, for a distribution, that its wave front
coincides with the generalized wave front of the embedded
distribution in $\mathcal{G}(\omega).$ Actually we prove a
stronger result,  the strong associativity of a distribution and a
generalized function is necessary and sufficient condition for the
equality of their analytic wave fronts. Section 4 is devoted to
the  comparison of the notions of S-analyticity in the sense of
\cite{ar2} and the analyticity of this paper.


\subsection{Definitions}
We recall the main definitions (\cite{co1},  \cite{gkos},
\cite{ober 001}, \cite{kuob}).
Let $\omega$ be an open set in $\mathbb{R}^{d}$
and $\mathcal{E}(\omega)$ be the  space of smooth functions with the
 sequence of seminorms
$\mu_\nu(\phi)=\sup\{|\phi^{(\alpha)}(x)|; \alpha\leq \nu,\; x\in
K_\nu\},\; \nu\in {\mathbb N}_0$, where $(K_\nu)_\nu$ is an
increasing  sequence of compact sets exhausting $\omega.$ Then the
set of moderate nets  $\mathcal{E}_{M}(\omega)$, respectively of
null nets  $\mathcal{N}(\omega),$  consists of nets $(f_{\varepsilon
})_{\varepsilon}\in \mathcal{E}(\omega)^{(0,1)}$ with the
properties
$$
(\forall n\in\mathbb{N})\;(\exists
a\in\mathbb{R})\;(\mu_{n}(f_{\varepsilon })=O(\varepsilon^{a})),
$$%
\[
\mbox{ respectively, } \;(\forall n\in\mathbb{N})\;(\forall
b\in\mathbb{R})\;(\mu_{n}(f_{\varepsilon})=O(\varepsilon^{b}))
\]
($O$ is the Landau symbol "big O").
 Both spaces are algebras and the latter is an ideal of  the former.
 Putting $v_n(r_\eps)=\sup\{a;\; \mu_n(r_\eps)=
O(\eps^a)\}$  and
 $e_n((r_\eps)_\eps,(s_\eps)_\eps)=\exp(-v_n(r_\eps-s_\eps)),\; n\in \N,$
we obtain  $(e_n)_n$,  a sequence of ultra-pseudometrics on
$\E_M(\omega)$ defining the ultra-metric topology (sharp topology) on
$\Ge(\omega).$
 The simplified Colombeau algebra
$\mathcal{G}(\omega)$ is defined as the quotient
$\mathcal{G}(\omega)=\mathcal{E}_M(\omega)/\mathcal{N}(\omega).$
This is also a differential algebra. If the nets
$(f_\varepsilon)_\varepsilon$ consist of constant functions on
$\omega $ (i.e. the seminorms $\mu_n$ reduce to the absolute value),
then one obtains
 the corresponding spaces of nets of complex  (or real) numbers $\mathcal{E}_M$ and
$\mathcal{N}_0.$ They are algebras, $\mathcal{N}_0$ is an ideal in
$\mathcal{E_M}$ and, as a quotient, one obtains the Colombeau
algebra of generalized complex numbers
$\bar{\mathbb{C}}=\mathcal{E}_M/\mathcal{N}_0$ (or
$\bar{\mathbb{R}}$). It is a ring, not a field. The sharp topology
in $\bar{\mathbb{C}}$ is defined as above. Note, a ball
$\tilde{B}(x_0,r)$, where
$x_0=[(x_{0,\eps})_\eps]\in\bar{\mathbb{C}}$ and $r>0,$ in the sharp
topology is given by $\tilde{B}(x_0,r)=\{[(x_\eps)_\eps] ; |x_{0,\eps}-x_\eps|=
O(\eps^{-\ln r})\}.$

Recall (\cite{gkos}), $\omega_M$ is the set of nets
$(x_\eps)_\eps\in\omega^{(0,1)}$ with the property
$|x_\eps|=O(\eps^{-a})$ for some $a>0$ and the relation of
equivalence is introduced by
$$(x_\eps)_\eps\sim (y_\eps)_\eps \mbox{ iff } |x_\eps-y_\eps|=O(\eps^a) \mbox{ for every } a>0.
$$
The corresponding classes of equivalence determine
$\tilde{\omega}$, and  $\tilde{\omega}_c$  denotes its subset
consisting of the classes of equivalence with the compactly
supported representatives; $(x_\eps)_\eps$ is compactly supported
by $K\subset \subset\omega$ if $x_\eps\in K$ for $\eps<\eps_0$.
The ultrametric in $\tilde{\omega}$ defined in the same  way as in
the one dimensional case for $\bar{\mathbb{C}}$. Note
$\tilde{\mathbb{C}}=\bar{\mathbb{C}}$
($\tilde{\mathbb{R}}=\bar{\mathbb{R}}$).

The embedding of Schwartz distributions in
$\mathcal{E}^{\prime}(\omega)$ is realized through the sheaf
homomorphism $ \mathcal{E}^{\prime}(\omega)\ni f\mapsto
[(f\ast\phi_{\varepsilon}|_{\omega})_\varepsilon]\in\mathcal{G}(\omega),
$ where the fixed  net of mollifiers
$(\phi_{\varepsilon})_{\varepsilon}$ is defined by
$\phi_{\varepsilon}=\varepsilon ^{-d}\phi(\cdot/\varepsilon),\;
\varepsilon<1,$ where $\phi\in\mathcal{S}(\mathbb{R}^{d})$
satisfies
$$\int
\phi(t)dt=1,\;\int
t^{m}\phi(t)dt=0,m\in\mathbb{N}_{0}^{n},|m|>0.$$
($t^{m}=t_{1}^{m_{1}}...t_{n}^{m_{n}}$ and $|m|=m_{1}+...+m_{n}$).
 $ \mathcal{E}^{\prime}(\omega)$ is embedded into
$\mathcal{G}_c(\Omega)$ of compactly supported generalized functions
and the
 sheaf homomorphism,  extended onto $\mathcal{D}^{\prime}$, gives the embedding of $\mathcal{D}^{\prime
}(\omega)$ into $\mathcal{G}(\omega).$

Generalized algebra  $\mathcal{G}_c(\omega)$  consists of elements
in $\mathcal{G}(\omega)$ which are compactly supported while
$\mathcal{G}^\infty(\omega)$ is defined in \cite{ober 001} as the
quotient  $\mathcal{E}^\infty_{M}(\omega)/\mathcal{N}(\omega),$
where $\mathcal{E}^\infty_{M}(\omega)$ consists of nets
$(f_{\varepsilon })_{\varepsilon}\in \mathcal{E}(\omega)^{(0,1)}$
with the property
$$
(\forall K\subset\subset \omega) (\exists a\in\mathbb{R}) (\forall
n\in\mathbb{N}_0^d) (|\sup_{x\in
K}f^{(n)}_{\varepsilon}(x)|=O(\varepsilon^{a})),
$$%
Note that $\mathcal{G}^\infty$ is a subsheaf of $\mathcal{G}.$

If $\Omega$ is an open set of $\mathbb{C}^{p}\equiv
\mathbb{R}^{2p}$ then a generalized function $f \in {\cal
G}(\Omega)$ is a holomorphic generalized function, if
$\bar{\partial}_{z_i}f=0$, $i=1,...,p,$ (\cite{col}) i.e. for any
representative $(f_\varepsilon)_\varepsilon$ of $f$,
$1/2(\partial_{x_i}+i\partial_{y_i})f_\varepsilon(x,y)\in
\mathcal{N}(\Omega)$, $i=1,...,p$.  In the sequel, we will write
$\bar{\partial}f=0$ for $ \bar{\partial}_{z_i}f=0$, $i=1,...,p.$
 It is proved in \cite{col}
that for every open bounded set of $W\subset\subset \Omega$ there
exists a representative $(f_\varepsilon)_\varepsilon$ of $f$ in
$W$ such that $f_\varepsilon$ is holomorphic in $W$ for every
$\varepsilon<1.$ By the sheaf property of ${\cal G}$, one obtains
that $\Omega\rightarrow {\cal G}_H(\Omega)$ is a sheaf. It is a
subsheaf of $\Ge^\infty.$

\section{Real analytic generalized functions}

\begin{definition}\label{defi, aaa}
Let $\omega$ denote an open set in ${\R^d}$ and $x_0$ be a point
of $\omega$. A function $g\in{\Ge}(\omega)$ is said to be real
analytic at $x_0$ if there exist an open ball
$B=B(x_0,r)\subset\omega$ containing $x_0$ and
$(g_\varep)_\varep\in {\E}_M(B)$ such that
$$
(i)\,\,\,\,\,\,\, f|B=[(g_\varep)_\varep]\,\,{\rm in}\,\,{\cal G}(B);
$$
$$
(ii)\,\,\,\,\,\,(\exists \eta>0)(\exists a>0)(\exists \varep_0\in
(0,1))
$$
$$
\sup_{x\in {B}}|g^{(\alpha)}_\varep(x)| \le
\eta^{|\alpha|+1}\alpha!\varep^{-a},\; 0<\varep<\varep_0,\;
\alpha\in\N^d.
$$

It is said that $g$ is real analytic in $\omega$ if $g$ is real
analytic at each point of $\omega$.
The space of all generalized functions which are real analytic in $\omega$ is denoted by 
$\mathcal{G}_A(\omega)$.

The analytic singular support, $\mbox{singsupp }_{ga} g,$ is the
complement of the set of points $x\in \omega$ where $g$ is real
analytic.
\end{definition}

\noindent It follows  from the definition that ${\cal G_A}$ is a subsheaf of ${\cal 
G}$.

Using Stirling's formula it is seen that condition (ii) in Definition \ref{defi, aaa} is 
equivalent to
$$
(iii)\,\,\,\,\,\,
(\exists \eta>0)(\exists a>0)(\exists \varep_0\in (0,1))
$$
$$
\sup_{x\in B}| g^{(\alpha)}_\varep(x)| \le
\eta(\eta|\alpha|)^{|\alpha|}\varep^{-a},\; 0<\varep<\varep_0,\;
\alpha\in\N_0^d.
$$
The use of Taylor expansion and condition $(ii)$ of Definition
\ref{defi, aaa} imply that $g_{\varep}, \eps<\varep_0$ admit
holomorphic extensions in a complex ball $B=\{z\in \mathbb{C}^d;
|z-x_0|<r \}$ which is independent of $\varep$. Consequently we
get a holomorphic extension $G$ of $[(g_\varep)_\varep]$ and then
$g|B\cap\mathbb{R}^d=G|B\cap\mathbb{R}^d$.

 The
following theorem gives a representative of a real analytic
generalized function in a neighborhood of
a point.
We will use a  multidimensional notation:

\noindent $\alpha!=\alpha_1!...\alpha_d!,\;
\eta^{\alpha}=\eta^{\alpha_1}...\eta^{\alpha_d}, \;\eta>0,\;
\alpha=(\alpha_1,...,\alpha_d)\in \mathbb{N}^d\;;$


\noindent $(iy)^j=(iy_1)^{j_1}...(iy_d)^{j_d}, i=\sqrt{-1},
y\in\mathbb{R}^d, j=(j_1,...,j_d)\in\mathbb{N}^d;$

\noindent In the one dimensional case (with $\zeta=\xi+i\eta,\; \xi,
\eta \in\mathbb{R}$),
 \\
$\partial f(\zeta)=1/2(\partial_\xi -i\partial_\eta) f(\xi,\eta)$,
 $\bar{\partial} f(\zeta)=1/2(\partial_\xi +i\partial_\eta) f(\xi,\eta).$

\begin{theorem}\label{prop, aaa}
Let $f\in{\cal{G}}(\omega)$  and
 $I\ssubset \omega$ be an open $d-$interval, $I=I_1\times...\times I_d$.
Let $\sigma$ denote a positive moderate function defined in $(0,1]$ ($\sigma(\eps)\leq p(1/\eps), \eps<1$
for some polynomial $p$)  such that
$\lim_{\varep\to 0}\frac{\sigma(\varep)}{\ln\varep}=-\infty$.
 Then there exists $\eta>0$ and $\varep_0\in (0,1)$ such
that a  net $(F_\varep)_\varep$ defined by
$$
F_\varep(x,y)=\sum_{j_1\le \sigma(\varep),..., j_d\le
\sigma(\varep)}f^{(j)}_{\varep}(x)(iy)^j/j!,\; (x,y)\in
I\times(-1/\eta,1/\eta)^d, \;\varep\in(0,\varep_0),
$$
represents a holomorphic generalized function $F=[(F_\eps)_\eps]$  in $I\times(-1/\eta,1/\eta)^d$
such that  $f=F(\cdot,0)$ in $I$ (and we say that $f$ is the boundary value of $F$).
\end{theorem}
\bproof
We will use the assumption that there exist $\eta>0,
\varep_0\in(0,1)$ and $a>0$ such that
$$\sup_{x\in I}|f^{(j)}_{\varep}(x)|\le \eta^{j+1}j!\varep^{-a},\; 0<\varep<\varep_0,\; j\in\N_0^d.$$

Take $ \varep\in (0,\varep_0)$ and $(x,y)\in
I\times(-1/\eta,1/\eta)^d.$ Later, as a new condition on $\eta$,
 we will increase $\eta$.
Let  $(f_\varep)_\varep$ be
a representative of $f$ in 
  $I$.

Let $m\in\N, m\leq d$, $j\in\N^d$ and $l\in\N^d$ such that $l_k=j_k$ for $k\ne m$ and $l_m=0$. Let $\hat{m}=(0,...,1,...0)$ 
(the $m$-th component equals $1$ and the others equal $0$). We
have
$$
2\overline{\partial}_{z_m}f^{(j)}_\varep(x)(iy)^j/j!=f^{(j+\hat{m})}_{\varep}(x)(iy)^j/j!
+if^{(j)}_{\varep}(x)ij_m(iy_m)^{{ j_m}-1}(iy)^{l}/j!
$$

$$
=(iy)^{l}/l!\left[f^{(j+\hat{m})}_{\varep}(x)(iy_m)^{j_m}/{{j_m}!}-
f^{(j)}_{\varep}(x)(iy_m)^{j_m-1}/({j_m}-1)!\right]$$
(if $j_m=0,$ then the second addend does not appear).

Let $E(\sigma(\varep)):=\max\{n\in\N; n\leq \sigma(\varep)\}$.
It follows
$$
2\sum_{j_k\le\sigma(\varep), k\leq d}\overline{\partial}_{z_m}f^{(j)}_\varep(x)(iy)^j/j!=
f^{(l+E(\sigma(\varep)+1)\hat{m})}_{\varep}(x)(iy_m)^{E(\sigma(\varep))}/E(\sigma(\varep)
)!.
$$
Thus,
$$
2\overline{\partial}_{z_m}F_\varep(x,y)=\sum_{l_k\le\sigma(\varep),k\ne m}(iy)^l/l!
f^{(l+E(\sigma(\varep)+1)\hat{m})}_{\varep}(x)(iy_m)^{E(\sigma(\varep))}/E(\sigma(\varep)
)!.
$$
Since
$$
\sup_{x\in I}\left|f_\varep^{(l+E(\sigma(\varep)+1)\hat{m})}(x)\right|\le 
\eta^{|l|+E(\sigma(\varep))+2}(l+E(\sigma(\varep)+1)\hat{m})!\varep^{-a},0<\varep<\varep_
0
$$
and $(l+E(\sigma(\varep)+1)\hat{m})!=l!E(\sigma(\varep)+1)!$, it follows that
$$
2\left|\overline{\partial}_{z_m}F_\varep(x,y)\right|\le E(\sigma(\varep)+1)
\sum_{l_k\le\sigma(\varep),k\ne m}
(|y^l|)\eta^{|l|+E(\sigma(\varep))+2}\varep^{-a}|y_m|^{E(\sigma(\varep))}.
$$
Now we enlarge $\eta$ and assume that $|y_k|<\rho/\eta, k=1,...,d, \;\rho<1$. We obtain
$$
2\left|\overline{\partial}_{z_m}F_\varep(x,y)\right|\le 
\eta^2\varep^{-a}\rho^{E(\sigma(\varep))}E(\sigma(\varep)+1)
\sum_{l_k\le\sigma(\varep),k\ne m}\rho^{|l|}.
$$
This gives
$$
2\left|\overline{\partial}_{z_m}F_\varep(x,y)\right|\le 
\eta^2(1-\rho)^{1-d}\varep^{-a}E(\sigma(\varep)+1)\rho^{E(\sigma(\varep))}.
$$
  With this and the
 condition  $\lim_{\varep\to 
0}\frac{\sigma(\varep)}{\ln\varep}=-\infty$
we have that for any $b>0$ there exists $\varep_b\in(0,1)$
such that $|\bar{\partial}_{z_m}F_\varep(x,y)|$ is bounded by $\varep^b$
for $\varep<\varep_b$. This concludes the proof.
\eproof
%

\begin{remark}\label{redkon}
Let $f\in\mathcal{G}_A(\omega).$ We will show in Theorem \ref{sss}
part b), that for every compactly supported generalized point
$x_0=[(x_{0,\eps})_\eps]$ ($x_0\in \tilde{\omega}_c$) and every
generalized point $x=[(x_\eps)_\eps]$ in a sharp neighborhood $V$
of  $x_0$   the series $F^N=[(F^N_\eps)_\eps], N\in \mathbb{N},$
converges to $f$ in the sense of sharp topology, where
$$
F^N_\varep(x_\eps)=\sum_{|j|=0}^N f^{(j)}_{\varep}(x_{0,\eps})(x_\eps-x_{0,\eps})^j/j! \;\;(\varep<\varep_0).
$$
\end{remark}

The existence of a global holomorphic representative of a holomorphic
generalized function in the one dimensional case (\cite{pili}), implies the following theorem.
\begin{theorem}
Let $f\in{\cal{G_A} }(\omega), \omega\subset \mathbb{R}$. Then $f$
admits a global real analytic representative, that is a
representative $(f_\varep)_\varep$ consisting of real analytic
functions such that $(i)$ and $(ii)$ in Definition 1 is satisfied
(with $f_\varep$ instead of $g_\varep$) for all $x_0\in\omega$ and
$B$ depending on $x_0$.
\end{theorem}
\bproof Let $x_0$ run over $\omega$. Choosing a representative which
satisfies (ii) in Definition \ref{defi, aaa}, we conclude that
there exists
 a ball  $B_{x_0}=B({x_0},r_{x_0})\subset\omega$,   a  complex
ball of the same radius $\tilde{B}({x_0},r_{x_0})=\tilde{B}_{x_0}$
 and $F_{x_0}\in \mathcal{G}_H({B}_{x_0})$ (defined with the same Taylor series
 as $f_\varepsilon, \varepsilon<\varepsilon_0$) such that 
$f|{B_{x_0}}=F_{x_0}|B_{x_0}$. As it is told, ${B}_{x_0}$ does
not depend on $\varep$. If $\tilde{B}_{x_0}\cap \tilde{B}_{x_1}\ne
\emptyset$, it follows that $F_{x_0}|\tilde{B}_{x_0}\cap
\tilde{B}_{x_1}=F_{x_1}|\tilde{B}_{x_0}\cap \tilde{B}_{x_1}$.
Actually, it is a consequence of Theorem \ref{ext} part b) given
below (or \cite{kelsca}, what is  part II quoted before Theorem
\ref{ext}).

Let $(f_{x_0,\varep})_\eps$ and $(f_{x_1,\varep})_\eps$ be representatives of $f$,
(with $B_{x_0}$ and $B_{x_1}$ as in Definition \ref{defi, aaa}
and let
 $x\in B_{x_0} \cap B_{x_1}$. Put $n_\varep(t)=f_{x_0,\varep}(t)-f_{x_1,\varep}(t),
t\in B_{x_0} \cap B_{x_1}, \varep \in (0,1).$
Since $(f_{x_0,\varep}-f_{x_1,\varep})_\eps$ is a net of real  analytic functions  in
$B_{x_0} \cap B_{x_1}$,   we have
$$ n_\varep(t) = \sum_{|j|=0}^\infty n^{(j)}_{\varep}(x)(t-x)^j/j!,\; t\in B(x,r), \eps \in (0,\eps_0),
$$
for some radius $r>0.$ Now, since 
 $\sup_{t\in B(x,r)}|n_\varep(t)|=O(\varep^a),$ for every $a>0,$
it follows that there exists $\eta>0$ such
that for every $a>0$ there exists $C>0$ and $\varep_0\in (0,1)$ such that
$$|n^{(j)}_{\varep}(x)|\leq C\eta^j j!\varep^a, \;j\in\mathbb{N},\; \varep<\varep_0.
$$
Consider $[(F_{x_0,\varep}-F_{x_1,\varep})_\eps]$ in the
intersection of complex balls $\tilde{B}_{x_0} \cap
\tilde{B}_{x_1}.$ It follows  that there exists an open complex
ball with the center at $x$ and of radius $1/\eta$, such that
$[(F_{x_0,\varep}-F_{x_1,\varep})_\eps]$ equals zero in this ball.
But we know (\cite{pili})
 that in this case $[(F_{x_0,\varep}-F_{x_1,\varep})_\eps]$
equals zero in $B_{x_0} \cap B_{x_1}$.

Let $\Omega=\cup_{x\in\omega}{B}_x$. Since $\Ge_H$ is a
sheaf, and  the family $F_x\in \Ge({B}_x), x\in \omega$ has
the property that on the intersections of the balls elements of the family
coincide, it follows that
there is $F\in{\cal G}(\Omega)$ such that for every $x\in \omega$ one has 
$F|{{B}_x}=F_x$. Since $F$ is
holomorphic, by \cite{pili} we have that  it admits a holomorphic representative 
$(\tilde{F}_\varep)_\varep$. Thus $(f_\varep)_\varep,$ where $f_\varep=\tilde{F}_\varep|\omega$, is a 
real analytic representative of $f$ which satisfies $(i)$ and
$(ii)$. \eproof

\begin{corollary}
a) Let $f\in\mathcal{G}_A(\omega),\; \omega\subset \mathbb{R}^d$. Then for every
open set $\omega_1\subset\subset\omega$ there exist an open set
$\Omega_1\subset\mathbb{C}^d$ such that $\Omega_1\cap\mathbb{R}^d=\omega_1$ and
 a net of holomorphic functions $(F_\eps)\eps\in\mathcal{E}_M(\Omega_1)$ such that
$F_\eps|\omega_1=f_\eps,\;\eps\in (0,1).$

b) If $d=1$ then the assertion in a) holds globally, with $\Omega$ instead of $\Omega_1.$
\end{corollary}

{\bf Hypothesis}. 1. The existence of a global holomorphic representative of $f\in\mathcal{G}_H(\Omega)$
depends on $\Omega\subset \mathbb{C}^d.$ Actually, we know such representations for appropriate domains
$\Omega\subset \mathbb{C}^n.$

2. If $f\in\mathcal{G}_A(\omega),\; \omega\subset \mathbb{R}^d$
then we conjuncture that it admits a global real analytic
representative.

\vspace{0,3cm}

It is proved in \cite{kelsca} that a  holomorphic generalized
function in an open set $\Omega\subset \mathbb{C}$
 is equal to zero if:

 {\it
 I It is equal to zero in every classical point of $\Omega;$

II It takes zero values at all points of a subset $S$ of a smooth
  curve $\gamma\subset \Omega$ with positive length measure, $meas{S}>0.$

III It takes zero values at all points of a subset $\mathcal{A}\subset\Omega$
of positive Lebesgue measure, $m(\mathcal{A})>0$.}
Clearly II $\Rightarrow$ I, and III $\Rightarrow$ I.

The proofs of (i), (ii) and (iii) of part a) of the next theorem is a
consequence of  quoted results and  the result of \cite{coga}
and \cite{pili}, where it is proved that a ($p-$dimensional)
holomorphic generalized function equals zero in $\Omega\subset \mathbb{C}^p$,
 if it equals zero in an open subset of this set.
 Actually, the basic assertion which will be used, is that a generalized function equals zero
 if and only if it equals zero at every compactly  supported generalized point \cite{kuob}.

 Although the first assertion in a) is the consequence of the second and third one, we will give  proofs
 of all assertions.
\begin{theorem} \label{ext}

 a) Let $\Omega$ be an open set of $\mathbb{C}^p, p>1$ and $f=[(f_\eps)_\eps]\in \Ge_H(\Omega).$

 (i)  Assume that
$f(x)=0$ for every  $x\in \Omega$ ($(f_\eps(x))_\eps\in\mathcal{N}(\Omega)$). Then $f\equiv 0.$

(ii) Assume that $S\subset\Gamma=\gamma_1\times\cdot\cdot\cdot\times\gamma_p\subset\Omega$, where $\gamma_i,i=1,...,p,$
are smooth curves in $\mathbb{C},$ has the property $meas(S)>0$, where $meas$ is the product of length measures,
and that $f(x)=0$ (in $\bar{\mathbb{C}}$) for every   $x\in S.$
Then $f\equiv 0.$

(iii) Assume that $\mathcal{A}$ is a Lebesgue measurable subset of
$\Omega,$ such that  its Lebesgue measure is positive
($m(\mathcal{A})>0)$ and that $f(x)=0$  for every   $x\in
\mathcal{A}.$ Then $f\equiv 0.$

b) Let $\omega$ be an open set of $\mathbb{R}^d, d\in\mathbb{N}$
and $f=[(f_\eps)_\eps]\in \Ge_A(\omega).$

If
$f(x)=0$ for every  $x\in \mathcal{A}$, where $\mathcal{A}\subset\omega$
is a Lebesgue measurable set and has a positive Lebesgues measure,
then $f\equiv 0.$
\end{theorem}
\bproof
 a)
 (i) We can assume that $\Omega$ is simply connected; otherwise we have to consider every
 component of $\Omega$ separately.
Denote by  $W=B_1\times...\times B_p$
 the  product of $p$   balls in $\mathbb{C}$ with radius $r>0$
  and center $t_i\in\mathbb{C}$, $i=1,...,p,$ such that
   $W\subset \subset \Omega.$ Using the sheaf property for $f\in\mathcal{G}_H(\Omega)$,
   it is enough to  prove that $f\equiv 0$ in $W$ and this is equivalent to show  that  $f(z)=0$
  for every compactly supported
   generalized point $z=(z_1,...,z_p)\in \tilde{W}_c$.

 Let us note that $z=(z_1,...,z_p)$ is a compactly supported generalized point
 of $\tilde{W}_c$ if and only if $z_i$ is a compactly supported generalized point of
 $\tilde{B_i},$
  for every $i=1,...,p.$

   Fix $(s_1,...,s_{p-1})\in B(t_1,r)\times...\times B(t_{p-1},r) $ and consider
   a holomorphic generalized function $$g_p(t):=f(s_1,...,s_{p-1},t),\; t\in B(t_p,r).$$ By assumption,
   it is equal to zero for every $t\in B(t_p,r)$. By \cite{kelsca} it is equal to zero in $B(t_p,r)$
   and  by \cite{kuob}, it is equal to zero in every compactly supported generalized
point $z_p\in(\tilde{B}(t_p,r))_c.$ Thus $f(s_1,...,s_{p-1},z_p)$
for every $(s_1,...,s_{p-1})\in B(t_1,r)\times...\times B(t_{p-1},r) $
and every compactly supported  $z_p$ in $\tilde{B}(t_p,r).$

Now we fix points $(l_1,...,l_{p-2})\in B(t_1,r)\times...\times B(t_{p-2},r) $, $z_p$ and consider
a holomorphic generalized function
$$g_{p-1}(t)=f(l_1,...,l_{p-2},t,z_p),\; t\in B(t_{p-1},r)$$
with the representative $g_{p-1,\eps}(t)=f(l_1,...,l_{p-2},t,z_{p,\eps}),\; t\in B(t_{p-1},r)$
By the same arguments, we have that $g_{p-1}(z_{p-1})=0$ for every
compactly supported generalized point $z_{p-1}\in (\tilde{B}(t_{p-1},r))_c$. Moreover, letting $z_p$ to run over
$(\tilde{B}(t_p,r))_c$ and  $z_{p-1}$ to run over
$(\tilde{B}(t_{p-1},r))_c,$
 we have that
$$f(l_1,...,l_{p-2},z_{p-1},z_p)=0,\; (z_{p-1},z_{p})\in (\widetilde{B(t_{p-1},r)\times B(t_p,r)})_c
$$
 for fixed but arbitrary $l_1,...,l_{p-2}.$ In this way we come to
$$g_1(t)=f(t,z_2,...,z_p), \;t\in B(t_1,r). $$  It is equal to zero at every compactly supported
 generalized point of $(\tilde{B}(t_1,r))_c$ and we conclude, that
$$f(z_1,...,z_p)=0,\; (z_1,...,z_p)\in (\widetilde{B(t_{1},r)\times...\times B(t_p,r)})_c.$$
 This finishes the proof of (i).

(ii)
The proof will be done by induction on the dimension $p.$ If $p=1$, this is II given above.

 Suppose that the assertion is true for the dimension $p-1.$ Let us write $z\in\mathbb{C}^p$
 as $z=(z^1,z^2),$ where $z^1\in\mathbb{C}^{p-1}$ and $z^2\in\mathbb{C}$.
 Let $s_1,...,s_p$ be subsets of $\gamma_1,...,\gamma_p,$ respectively, with positive one dimensional
 length measure such that $(S^1\times S^2)=((s_1\times...\times s_{p-1})\times s_p)\subset S$ and let
 $B_1,...,B_p$ be one dimensional balls such that $s_i\subset B_i, i=1,...,p.$ and the product of
 balls is contained in $\Omega.$
 Fix $z^1\in S^1.$ Then $g_{z^1}(t)=f(z^1,t), t\in B_1$ is a holomorphic generalized function
 such that $g_{z^1}(t)=0$ for every $t\in s_1$ and by assertion II, $g_{z^1}\equiv 0.$  Thus,
 it holds for every $z^1\in S^1$ and every compactly supported  generalized point $\tilde{t}
 \in B_p.$ Now fix $\tilde{t}$ a compactly supported point of $B_p$ and consider $f_{\tilde{t}}(z^1),
 z^1\in B_1\times...\times B_{p-1}.$ By induction,
 we have that it is equal to zero at every compactly supported point of
$ B_1\times...\times B_{p-1}$.  This completes the proof that $f$ is equal to zero at every
compactly supported point of $B_1\times...\times B_{p}$ and we have $f\equiv 0$
in $\Omega$.

(iii) Consider a polydisc $\Delta=D_1\times...\times D_p$ such
that $m(\mathcal{A}\cap \Delta)>0.$ Now using assertion III, as
above, we have that $f$ is equal to zero in $\Delta,$ hence in
$\Omega$.

b)
 There exists a box $I=I_1\times...\times I_d\subset \omega$
such that $m(\mathcal{A}\cap I)>0.$ Note $I_i$ are curves in $\mathbb{C}.$ Consider now
the holomorphic extension of $f$ in $\Omega\subset\mathbb{C}^d$, where $\Omega\cap \mathbb{R}^d=I.$
Thus $g(x)=0 $ for every $x\in \mathcal{A}\cap I$ and this implies $g\equiv 0$ in $\Omega$ (part (ii) of a)), hence
$f=0$ in $I$.
 This completes the proof of part b).
 \eproof

\begin{theorem} Let $f\in\mathcal{E}'(\omega)$ and $f_\varepsilon=f*\phi_{\varepsilon}, \varepsilon\in (0,1)$
be its regularization by a net $\phi_\varepsilon=\frac{1}{\varepsilon}\phi(\cdot/\varepsilon), \varepsilon<1,$
where $(\phi_\varepsilon)_\varepsilon$ is a net of mollifiers. Then
$$\mbox{ singsupp }_a f = \mbox{ singsupp }_{ga} [(f_\varepsilon)_\varepsilon],
$$
where on the left hand side is the analytic singular support of the distribution $f$.
\end{theorem}
The proof is the consequence of the proof of  Theorems \ref{444} and \ref{555}.

\begin{example}
\em  We will construct a generalized function by a net of real
analytic functions  which is not  a real analytic generalized
function.

Let $f_\varepsilon(x)=\frac{t}{\cosh(t/\varepsilon)}, \; t\in
\mathbb{R},\; \varepsilon\in (0,1).$ This net determines a
generalized function $f\in\mathcal{G}(\mathbb{R})$ such that for
every $t\in \mathbb{R},$ $f(t)=0$ in the sense of generalized real
numbers $\bar{\mathbb{R}}$ but $f\neq 0$ in
$\mathcal{G}(\mathbb{R}).$ Moreover, $f^{(\alpha)}(t)=0$ in
$\overline{\mathbb{R}}$, for every $t\neq 0$.

 For every fixed $\varepsilon,$
$f_\varepsilon$ is real analytic, but there does not exist a
common open set $V$ around  zero
such
that $f_\varepsilon$ are analytic in $V$ for
$\varepsilon<\varepsilon_0.$ More precisely, this generalized function is not real analytic at zero.

Note, by $f_\varepsilon(z)=\frac{z}{\cosh(z/\varepsilon)},\;
\varepsilon\in (0,1)$ is defined a moderate net in
$\mathbb{C}\setminus\{z; \Im z=0\}$ and $[(f_\varepsilon(z))]=0$
in $\mathbb{C}\setminus\{z; \Im z=0\}$.
\end{example}

\section{Analytic wave front set}
We will use the notation $\hat{\phi}$ for the Fourier transformation $\mathcal{F}(\phi)(\xi)
=\int_{\mathbb{R}^d}e^{-ix\xi}\phi(x)dx,\; \xi\in\mathbb{R}^d, \; \phi \in
L^1(\mathbb{R}^d)$.

\begin{definition}\label{d3}
Let $\omega$ be an open set in ${\R^d}$. A sequence $(u_n)_n$ in ${\cal G}_c(\omega)$ will 
be called bounded if there exists  a sequence of representatives
$(u_{n,\eps})_\eps$, $n\in\mathbb{N},$ such that
\begin{equation} \label{nnn}
(\exists K\ssubset\omega)(\exists m\in\R)(\exists b>0)
\end{equation}
$${\supp}(u_{n,\varep})\ssubset K, \varep\in(0,1] \;\;\mbox{ and}
\;\;
\sup_{\xi\in\R^d}(1+|\xi|)^{-m}|\widehat{u_{n,\varep}}(\xi)|=O(\varep^{-b}).
$$
\end{definition}

The following lemma, related to a special sequence of smooth functions
$(\kappa_n)_n$
 is a consequence of \cite{hor}, Theorem 1.4.2.
We will need the next two  interpretations of such a sequence and we refer to  \cite{hor}
for the proof.
\begin{lemma}\label{lemm, aaa}

 (i) Let
$K$ be a compact set 
of $\R^d$, $r>0$ and  $K_r=\{x:d(x,K)\le r\}.$ There exist
$C>0$ 
and a sequence of smooth
functions $(\kappa_n)_n$ , such that:
\\
(a) $\kappa_n=1$ on $K$ and ${\supp}(\kappa_n){\ssubset} K_r$ for
all $n\in\N$,
\\
(b)  $\sup_{x\in\mathbb{R}^d, |\alpha|\le
n}|\kappa^{(\alpha)}_{n}(x)|\le C(C n/r)^{|\alpha|}$, $n\in\N$.
\\
(ii) For every $x\in \omega $ there exist open sets $V, W$, $x\in
W\subset\subset V\subset \subset \omega,$  and a sequence of
functions $(\kappa_n)_n$ in $\mathcal{D}(V)$ with the property
$\kappa_n\equiv 1$ in $W$, $n\in \mathbb{N},$ and that for every
$\alpha\in\mathbb{N}^d$ there exists $C_\alpha>0,$ such that
$$|\kappa_n^{(\alpha+\beta)}(x)| \leq C_\alpha ( C_\alpha n)^\beta,\; |\beta|\leq n,\; n\in \mathbb{N}.
$$
\end{lemma}
We have the following characterization of analytic generalized functions.
\begin{theorem}\label{th, aab}
Let  $x_0\in\omega\subset {\R^d}$. An $f\in{\cal 
G}(\omega)$
is an analytic generalized function at $x_0$ if and only if there exist a bounded sequence $(u_n)_n$ in ${\cal 
G}_c(\omega)$
and $B(x_0, r)\subset \omega$ such that the following two 
conditions hold:
\\
(i) $u_n=f$ in $B(x_0, r)$ for all $n\in\N$;
\\
(ii)There exists a sequence of  representatives
$(u_{n,\varep})_\varep\in {\cal E}_{M,c}(\omega)$,
$n\in\mathbb{N},$ which satisfies
$$(\exists C>0)(\exists a>0)(\exists \varep_0\in (0,1])
$$
$$
|\widehat{u_{n,\varep}}(\xi)|\le (Cn)^n(1+|\xi|)^{-n}\varep^{-a},
\xi\in\R^d, 0<\varep<\varep_0, n\in\N.
$$
\end{theorem}
\bproof
Assume that $f$ is analytic at $x_0$. With the notation of Definition \ref{defi, aaa}, let 
$B=B(x_0,\lambda)$
for some $\lambda>0$. Using condition (iii) of Definition \ref{defi, aaa},
 there exist $\eta_1>0,a>0,\varep_0\in(0,1)$ 
such that
\begin{equation} \label{equa, aaa}
\sup_{x\in {B}}| f^{(\alpha)}_\varep(x)| \le
\eta_1(\eta_1|\alpha|)^{|\alpha|}\varep^{-a},0<\varep<\varep_0,\alpha\in\N^d.
\end{equation}
Choose $r=\lambda/3$ and take the sequence $(\kappa_n)_n$ of Lemma \ref{lemm, aaa} (i), associated 
to
$K=B(x_0,r)$. Set $u_{n,\varep}=\kappa_n f_\varep, n\in\N,\varep\in(0,1)$. Then $u_{n,\varep}, 
n\in\N,\varep\in(0,1)$,
define a sequence $(u_n)_n$ in ${\cal G}_c(\omega)$ such that ${\supp}(u_{n,\varep})\subset 
B(x_0,2r)$ for $n\in\N, \varep\in(0,1]$. We have ($\xi
\in\mathbb{R}^d$)
\begin{equation} \label{equa, aab}
\widehat{u_{n,\varep}}(\xi)=\int_{\R^d}
e^{-ix.\xi}u_{n,\varep}(x)dx= \int_{B(x_0,2r)}
e^{-ix.\xi}\kappa_{n}(x)g_\varep(x)dx.
\end{equation}
 It follows that for $\eps<\eps_0$ and  $\xi\in\mathbb{R}^d$
$$
|\widehat{u_{n,\varep}}(\xi)|\le \mbox{ mes }(B(x_0,2r))\sup_{x\in
B(x_0,2r)} |\kappa_{n}(x)g_\varep(x)|,\;n\in\mathbb{N}.
$$
Condition $(b)$ of Lemma \ref{lemm, aaa} (i) and (\ref{equa, aaa}) with $\alpha=0$ 
give
a constant $c>0$ such that
$$
|\widehat{u_{n,\varep}}(\xi)|\le c\varep^{-a}, \;\xi\in\R^d,
0<\varep<\varep_0, n\in\mathbb{N}.
$$
Thus $(u_n)_n$ is a bounded sequence in ${\cal G}_c(\omega)$ and obviously, $u_n=f, n\in \N$ in 
$B(x_0,r)$.
\\
For a given $n\in\N$ and $|\alpha|\le n$, we have
${u^{(\alpha)}_{n,\varep}}=
\sum_{\beta\le\alpha}(_\beta ^\alpha)
f^{(\beta)}_\varep\kappa^{(\alpha-\beta)}_n$.
\\
Using again $(b)$ of Lemma \ref{lemm, aaa} (i) and  (\ref{equa, aaa}) we have that there exists  
 $\eta>0$
such that
\begin{equation}
\sup_{x\in\mathbb{R}^d}|{u^{(\alpha)}_{n,\varep}}(x)|\le\eta(\eta n/r)^{|\alpha|}\varep^{-a}, 
n\in\N, |\alpha|\le n,
0<\varep<\varep_0.
\end{equation}
Since $\xi^\alpha\widehat{u_{n,\varep}}(\xi)=\int_{B(x_0,2r)} e^{-ix.\xi}
u^{(\alpha)}_{n,\varep}(x)dx$, we have, for $\xi\in\mathbb{R}^d,$
$$
|\xi^\alpha||\widehat{u_{n,\varep}}(\xi)|\le{\mes}{(B(x_0,2r))}\eta(\eta 
n/r)^{|\alpha|}\varep^{-a},
|\alpha|\le n,0<\varep<\varep_0.
$$
This implies, with  $|\alpha|=n,$
$$
|\xi|^n|\widehat{u_{n,\varep}}(\xi)|\le{\mes}{(B(x_0,2r))}\eta(\eta n/r)^{n}\varep^{-a}, 
\xi\in\R^d, n\in\N, 0<\varep<\varep_0.
$$
It follows that there exists  $C>0$ such that
$$
(1+|\xi|)^n|\widehat{u_{n,\varep}}(\xi)|\le (Cn)^{n}\varep^{-a} ,\xi\in\R^d, 
n\in\N,0<\varep<\varep_0.
$$
Conversely, assume that the sequence $(u_n)_n$ in Theorem 5 exists. Then
for $\eps<\eps_0$ and $x\in B(x_0,r)$ we 
have
$$ u^{(\alpha)}_{n,\varep}(x)=(2\pi)^{-d}\int_{\R^d} 
e^{ix.\xi}\xi^\alpha\widehat{u_{n,\varep}}(\xi)d\xi,\;
 |\alpha|\le n.$$
\\
For a given $\alpha$  we choose $n=|\alpha|+d+1$ and get
$$
| u^{(\alpha)}_{n,\varep}(x)|\le (2\pi)^{-d}(Cn)^n\varep^{-a}\int_{\R^d} 
|\xi^\alpha|(1+|\xi|)^{-n}d\xi, x\in B(x_0,r), \eps<\eps_0.$$
We have $\int_{\R^d} |\xi^\alpha|(1+|\xi|)^{-n}d\xi\le \int_{\R^d} 
(1+|\xi|)^{-d-1}d\xi<\infty$.
It follows that there exists a constant $C_1$ such that (with $n=|\alpha|+d+1$)
$$
|u^{(\alpha)}_{n,\varep}(x)|\le
(C_1|\alpha|)^{|\alpha|}\varep^{-a}, x\in B(x_0,r), \eps<\eps_0.
$$
Now, since $f_\varep={u_{n,\varep}}_{|{\omega_0}},$ the proof of the proposition
is completed.
\eproof

Theorem \ref{th, aab} leads to the following definition.

\begin{definition} \label{d4}
Let $f\in{\cal G}(\omega)$. Then, the analytic wave front
set of $f$, $\wf_{ga}(f)$, is the 
closed subset of $\omega\times\R^n\setminus \{0\}$
whose complement consists of points $(x_0, \xi_0)\in \omega\times\R^n\setminus \{0\}$ 
satisfying the following conditions:

There exist $B(x_0,r)$, an open conic neighborhood $\Gamma$ of 
$\xi_0$ and a bounded sequence $(u_n)_n$ in ${\cal G}_c(\omega)$
such that $u_n=f, $ in $B(x_0,r), n\in\mathbb{N},$ with
a sequence of representatives $(u_{n,\varep})_\varep, n\in \N,$
with the property
$$(\exists a>0)(\exists \varep_0\in(0,1])(\exists C>0)$$
$$|\widehat{u_{n,\varep}}(\xi)|\le (Cn)^n(1+|\xi|)^{-n}\varep^{-a},\; \xi\in\Gamma, 
0<\varep<\varep_0, n\in\N.$$
\end{definition}

By the use of sequence $(\kappa_n)_n$ of Lemma \ref{lemm, aaa}, in the same way as in \cite{hor}, Lemma
8.4.4, one can show that the next definition is equivalent to the
previous one.

\begin{definition}\label{d5}
Let
$f\in{\cal G}(\omega)$ and $(x_0,\xi_0)\in \omega\times \mathbb{R}^d\setminus\{0\}.$
 Then, it  is said that $f$ is $g-$microanalytic at this point if:
 \begin{itemize}
 \item[]
There exist
 neighborhoods $W$ and $V$ and a   sequence $(\kappa_n)_n$ as in
Lemma \ref{lemm, aaa}, (ii);

\item[]
There exists a cone $\Gamma_{\xi_0}$ around $\xi_0;$
\item[]
There exist $a>0$, $C>0$ and $\varep_0,$
\end{itemize}
such that
\begin{equation} \label{vf}
|\widehat{\kappa_nf_{\varepsilon}}(\xi)|\le C\varepsilon^{-a} (\frac{C n}{1+|\xi|} )^n, \; n\in\mathbb{N},\;
\xi\in \Gamma,\;  \varep\in(0,\varep_0).
\end{equation}
Then, the analytic wave front set, $\mbox{WF }_{ga} f$ is the
complement of the set of points where $f $  is $g$-microanalytic.
\end{definition}

Recall (\cite{nps}), $T\in\mathcal{D}'(\omega)$ and $f=[(f_\varepsilon)_\varepsilon]
\in\mathcal{G}(\omega)$ are strongly associated if  for every $K\subset\subset \Omega$ there exists $a>0$ such that
$|\langle T-f_\varepsilon,\theta\rangle|= O(\varepsilon^a), \varepsilon\rightarrow 0,$ for every
$\theta\in \mathcal{D}(\omega)$ with supp$ \theta\subset K.$
If $|\langle T-f_\varepsilon,\theta\rangle|\rightarrow 0, \varepsilon\rightarrow 0$,
 then $T\in\mathcal{D}'(\omega)$ and $f=[(f_\varepsilon)_\varepsilon]$ are associated (\cite{col}).

Also recall that  a distribution $f\in\mathcal{D}'(\omega)$ is
microanalytic  at $(x_0,\xi_0)\in \omega\times
\mathbb{R}^d\setminus\{0\}$ if (\ref{vf}) holds for
$f=f_\varepsilon$ on the left-hand side and $a=0$ on the
right-hand side.

The next example serves as a good motivation for  Theorem \ref{str}.
\begin{example}
\em

Let $\hat{\psi}  \in \mathcal{D}(\mathbb{R})$ equal one in the unit ball and let
 $\psi$ be its inverse Fourier transform.
Define $g_\varepsilon(x) = |\log \varepsilon|
\psi(x|\log\varepsilon|), x\in\mathbb{R},$ $\varepsilon<1.$ This
net determines  a real analytic generalized function in
$\mathbb{R}$. To prove this, one has to use the inequality
$|\varepsilon (\log \varepsilon)^n| \leq e^{-n} n^n,
\varepsilon<1,$ in order to estimate  $|g^{(n)}_\varepsilon|$.

This analytic generalized is  associated to the $\delta$
distribution, thus not to  an analytic distribution.

Moreover, if we start with a compactly supported distribution $T$ and consider
  $T_\varepsilon(x) = T*|\log \varepsilon| \psi(x|\log\varepsilon|), x\in\mathbb{R},$ $\varepsilon<1,$
  we obtain a real analytic generalized function $[(T_\varepsilon)_\varepsilon].$
  It is associated to $T$ but not strongly associated to $T.$
 \end{example}
\begin{theorem} \label{str}
Let $f=[(f_\varepsilon)_\varepsilon]\in\mathcal{G}(\omega)$ be  $g-$microanalytic  at
$(x_0,\xi_0)\in \omega\times \mathbb{R}^d\setminus\{0\}$ and let
 $T\in \mathcal{D}'(V)$ be strongly associated to
$f_{|V},$  where $V$ is an open neighborhood of $x_0$, $V\subset \omega$.
  Then $T$ is  microanalytic at $(x_0,\xi_0).$
\end{theorem}
\bproof
We will assume that $\varepsilon_0=1$ in (\ref{vf}). Moreover, put $f_\varepsilon=0, \varepsilon\geq 1.$
Then (\ref{vf}) holds for $\varepsilon>0.$ This simple remark will be important in the proof.

Assume that $\kappa_n$ equals one in  $W$ and that it is supported
by $V, n\in \N$ (cf. Lemma \ref{lemm, aaa} (ii)). Since $(\kappa_n
T)_n$ is a bounded sequence in $\mathcal{E}'(V)$ we know that it
is bounded in  the Banach space $(C^k(\overline{V}_1))'$, where
$\omega \supset V_1\supset V,$ for $k\in \mathbb{N}$ which is
equal to the order of $T$ (see Lemma 8.4.4 in \cite{hor}). Thus we
have that there exists $M>0$ such that
\begin{equation}
\label{ogr} |\langle T,\rho \rangle|\leq M\sup_{|\alpha|\leq k,
x\in V_1}|\rho^{(\alpha)}(x)|.
\end{equation}
By the strong association, there exists $a>0$ such that $|\langle
T-f_\varepsilon,\theta\rangle|= O(\varepsilon^a),
\varepsilon\rightarrow 0$ for every $\theta\in \mathcal{D}(V_1)$.
Let $a=b+c, b,c>0.$ We have that
$$T_{b,\varepsilon}=\varepsilon^{-b}(T-f_\varepsilon) \rightarrow 0 \mbox{ in } \mathcal{D}'(V_1)
\mbox{ as } \varepsilon\rightarrow 0  $$ and thus
$T_{b,\varepsilon}\rightarrow 0$ in $(C^k(\overline{V}_1))'$,  as
$\varepsilon\rightarrow 0,$ for some $k\in \mathbb{N}$. We assume
that this is the same $k$ as in (\ref{ogr}). By the Banach-Steinhaus
theorem, $\{T_{b,\varepsilon}; \varepsilon \in (0,1)\}$
is bounded in $(C^k(\overline{V}_1))'$ and there exists $M>0$ such
that for every $\rho\in C^k(\overline{V}_1)$
\begin{equation}\label{ogr1}
|\langle T-f_\varepsilon,\rho\rangle|\leq M\varepsilon^b
\sup_{|\alpha|\leq k, x\in V_1}|\rho^{(\alpha)}(x)|,\varepsilon<1
\end{equation}
Let us write (\ref{ogr}) in the form
\begin{equation}\label{ogr2}
|\langle T,\rho\rangle|\leq M\varepsilon^b\sup_{|\alpha|\leq k, x\in V_1}|\rho^{(\alpha)}(x)|,
\varepsilon\geq 1.
\end{equation}
Since $|\kappa_n^{(\alpha)}|\leq C(Cn)^{|\alpha|}, |\alpha|<n$  (Lemma \ref{lemm, aaa} (ii)),
we have, for $\varepsilon<1,$
$$|\widehat{\kappa_n(T-f_\varepsilon)}(\xi)|
\leq M\varepsilon^b\sup_{|\alpha|\leq k, x\in V_1}|(\kappa_n(x)e^{-ix\xi})^{(\alpha)}|
$$
and
$$|\widehat{\kappa_nT}(\xi)|\leq |\widehat{\kappa_nf_\varepsilon}(\xi)|+
 C\varepsilon^b (Cn)^k(1+|\xi|)^{k},\;\xi\in\mathbb{R}^d.
$$
 By (\ref{vf}), we have
$$|\widehat{\kappa_nT}(\xi)|\leq C(Cn)^n(1+|\xi|)^{-n}\varepsilon^{-a}+
 C\varepsilon^b (Cn)^k(1+|\xi|)^{k},\;\xi\in\Gamma, \varepsilon<1.
$$
Now put
\begin{equation}
\label{epsc}
\varepsilon=C^{-1/b}(Cn)^{-k/b}p^{p/b}(1+|\xi|)^{-(p+k)/b}, \;\xi\in \Gamma, p\in \mathbb{N}, n\in\mathbb{N}.
\end{equation}
We will consider separately two possible cases, when such defined $\varepsilon$ satisfies
$\varepsilon<1$ and $\varepsilon\geq 1.$
If $\varepsilon<1,$ we continue,
$$|\widehat{\kappa_nT}(\xi)|\leq
\frac{(Cn)^n(1+|\xi|)^{-n+(p+k)a/b}A^{a/b}(Cn)^{ka/b}}{p^{p+ap/b}}+C(Cp)^p(1+|\xi|)^{-p}.
$$
Put $n=p(a/b+1)+ka/b$.  Then  $n^n/p^{p(1+a/b)}\leq r^p$ for suitable $r>0$ and
with the new constants $C$ and $p$ so that  $\varepsilon<1 $
in (\ref{epsc}), we have
$$|\widehat{\kappa_{p(1+a/b)+ka/b}T}(\xi)|\leq
C(Cp)^p(1+|\xi|)^p, \xi \in \Gamma.$$
Thus,
 putting $\psi_p=\kappa_{p(1+a/b)+ka/b}$ for those $p$ for which $\varepsilon<1$ in (\ref{epsc}),
 we obtain
\begin{equation}
\label{kr}
|\widehat{\psi_{p}T}(\xi)|\leq
C(Cp)^p(1+|\xi|)^p, \xi \in \Gamma.
\end{equation}
Now we will consider those $p$ and $n=p(a/b+1)+ka/b$ for which  $\varepsilon\geq 1$ in (\ref{epsc}).
Now consider (\ref{ogr2}). We have immediately that (\ref{kr}) holds.
This completes the proof.
\eproof

The next theorem corresponds to H\" ormander's Theorem 8.4.5 in \cite{hor}, and it
is already proved in \cite{mar}. In fact in \cite{mar} is considered the wave front
which corresponds to a  general sequence $(L_k);$ in our case $L_k=k$.
We just quote it here.
\begin{theorem}\label{444}
Let
$f\in{\cal G}(\Omega)$. Then $pr_1 \mbox{WF }_{ga} f=\mbox{singsupp }_{ga} f.$
\end{theorem}
Now we will compare the wave front set for an $f$ in $\mathcal{E}'(\omega)$ and the wave front of the corresponding
generalized function.

\begin{theorem}\label{555}
Let
$f\in {\cal E}'(\omega)$ and
$f_\varep =f*\phi_\varep|_\omega,\; \varep\in (0,1).$
(It is  a representative of embedded distribution; $(\phi_\eps)_\eps$  is a net of
mollifiers  described in the introduction).
Put  $F=[(f_\varep)_\varep].$
Then
$$\mbox{WF }_{a}f = \mbox{WF }_{ga} F,$$
where on the left side is the analytic wave front set of  the distribution $f$.
\end{theorem}
\bproof
Since $f$ and $F$ are strongly associated, previous theorem implies that $WF_{ga}F\supset WF_a f.$
We will show the opposite inclusion by showing that there exists a representative $(g_\varepsilon)_\varepsilon$
of $F$ with the property: If $f$ is microanalytic at $(x_0,\xi_0)$, then
\begin{equation}
\label{ho}
|\widehat{\kappa_ng_{\varepsilon}}(\xi)|\le C\varepsilon^{-a} (\frac{C n}{1+|\xi|} )^n, \;
\xi\in \Gamma,\;  \varep\in(0,\varep_0),\; n\in\mathbb{N}.
\end{equation}
Let $h$ be a nonnegative function of $\mathcal{D}(\mathbb{R}^d)$ such that it is equal one in the ball $B(0,1)$
and which is supported by $B(0,2).$
Let
$$ \psi_\varepsilon(x)=h(\varepsilon^{-1/2}x)\phi_\varepsilon(x), x\in\mathbb{R}^d, \varepsilon \in (0,1).
$$
We will prove later the following lemma.
\begin{lemma}
\label{dod} For every $k,\alpha\in\mathbb{N}^d$ and every $a>0$
$$\sup_{x\in\mathbb{R}^d}|x^k(\phi^{(\alpha)}_\varepsilon(x)-
\psi^{(\alpha)}_\varepsilon(x))| =O(\varepsilon^a).$$
\end{lemma}
Thus, by this lemma, $(g_\varepsilon)_\varepsilon=(f*\psi_\varepsilon)_\varepsilon$ is also a representative of
$f$. We will show that (\ref{ho})
holds for this net.

Let $(\kappa_n)_n$ and $(\tilde{\kappa}_n)_n$ be sequences with the properties of Lemma \ref{lemm, aaa} (ii).
so that $\tilde{\kappa}_{n}(x)=1, x\in $supp$\kappa_{n+d+1}, n\in \mathbb{N}.$
By the estimate of Lemma \ref{lemm, aaa}, we have that there exists $M>0$ which does not depend on $n$
such that
$$\int_{\mathbb{R}^d}|\widehat{\kappa_{n+d+1}}(t)|(1+|t|)^n dt<M, n\in \mathbb{N}.
$$

There exists $\varepsilon_0\in (0,1)$ such that
$$\kappa_{n+d+1}g_\varepsilon=\kappa_{n+d+1}(f\tilde{\kappa}_n*\psi_\varepsilon), \varepsilon \in(0,\varepsilon_0), n\in \mathbb{N}.$$
Thus, by Peetre's inequality, for $\varepsilon \in
(0,\varepsilon_0),$
$$|\widehat{\kappa_{n+d+1}g_\varepsilon}(\xi)|\leq
\int_{\mathbb{R}^d}|\widehat{\kappa_{n+d+1}}(\xi-t)\widehat{f\tilde{\kappa}_n}(t)
\widehat{\psi}(\varepsilon t)|dt
$$
$$\leq M C(Cn)^n\int_{\mathbb{R}^d}|\widehat{\kappa_{n+d+1}}(t)|(1+|\xi-t|)^{-n}dt\leq \frac{C(Cn)^n}{(1+|\xi|)^{-n}}
\int_{\mathbb{R}^d}|\widehat{\kappa_{n+d+1}}(t)|(1+|t|)^ndt.
$$
This completes the proof of the theorem.
\eproof
{\bf Proof of Lemma \ref{dod}}.
Put $s_\varepsilon=\phi_\varepsilon-\psi_\varepsilon,\varepsilon<\varepsilon_0.$ Note that $s_\varepsilon(x)=0 $
if $|x|<\varepsilon^{1/2}, \varepsilon<\varepsilon_0.$

Let  $k, \alpha\in\mathbb{N}_0^d$ and
$|x|\geq  \varepsilon^{1/2},
\varepsilon<\varepsilon_0.$
We have
$$|x^k s^{(\alpha)}_\varepsilon(x)|\leq |x|^{|k|}\sum_{\beta+\gamma=\alpha}| (1- h(\varepsilon^{-1/2} x))^{(\beta)}
|\phi^{(\gamma)}_\varepsilon(x)|
$$
$$\leq C |x|^{|k|}\sum_{\beta+\gamma=\alpha}c_{\beta,\gamma}\varepsilon^{-\beta/2}\varepsilon^{-|\gamma|+d}
|\phi^{(\gamma)}(x/\varepsilon)|
$$
$$\leq
C |x|^{|k|}\varepsilon^{-|\alpha|-d}
\sup_{\gamma\leq\alpha}|\phi^{(\gamma)}(x/\varepsilon)|.
$$
(recall again, $\phi$ defines a net of mollifiers $(\phi_\varepsilon)_\varepsilon$.)
Put $u=x/\varepsilon.$  The last expression becomes

$$\varepsilon^{|k|-(d+|\alpha|)}\sup_{\gamma\leq\alpha}|u|^{|k|}|\phi^{(\gamma)}(u)|, |u|\geq
\varepsilon^{-1/2}, \varepsilon<\varepsilon_0.
$$
Since $\phi\in\mathcal{S}(\mathbb{R}^d),$  for every $m,\alpha\in\mathbb{N}^d,$
there exists $C_{m\alpha}>0$
such that
$$|u^{2|m|+|k|}\phi^{(\alpha)}(u)|\leq C_{m,\alpha}, u\in \mathbb{R}^d,$$
which implies
$$|u|^{|k|}|\phi^{(\alpha)}(u)|\leq C\varepsilon^{-2|m|}, |u|\geq \varepsilon^{-1/2},
\varepsilon<\varepsilon_0.
$$
Thus we have
$$|x^k s^{(\alpha)}_\varepsilon(x)|\leq \varepsilon^{|k|-(d+|\alpha|)+|m|}, |x|>\varepsilon^{1/2}
$$
and choosing enough large $m$, we have that the right hand side is equal to $O(\varepsilon^a)$ for every
$a>0$. This finishes the proof of lemma.


\section{S-analytic generalized functions }

S-analytic generalized functions are
introduced in \cite{ar2} through the Taylor expansions in the sharp topology of
$\tilde{\Omega}\subset \tilde{\mathbb{C}}$. Actually, in \cite{ar2} such generalized functions are
 called  "analytic" but since this definition leads to
a quite different algebra of generalized functions, we will call
them S-analytic, in order to avoid the confusion.

 Here we will consider the corresponding definitions for
$\tilde{\Omega}\subset \tilde{\mathbb{R}}^d$
(if $d\equiv 2p,$ we use notation $\mathbb{C}^p=\mathbb{R}^d$) and we will distinguish the cases
$\tilde{\Omega}\subset \tilde{\mathbb{C}}^p$ and  $\tilde{\omega}\subset \tilde{\mathbb{R}}^d$.

 We reformulate the definition of S-analyticity of \cite{ar2}.
\begin{definition} \label{san} (\cite{ar2})
Let $\omega\subset \mathbb{R}^d$ (resp., $\Omega\subset \mathbb{C}^p$) be  open and $\tilde{s}_0\in
\tilde{\omega}_c$, (resp., $\tilde{s}_0\in
\tilde{\Omega}_c$). Then   $f: \tilde{\omega} \rightarrow
\tilde{\mathbb{C}}$ is S-real analytic at $\tilde{s}_0$ (resp.,  S-holomorphic at $\tilde{s}_0$) if there exist a
sequence $a_n \in\tilde{\mathbb{C}}, n\in\mathbb{N}_0^d$  and a series of the
form $\sum_{|n|=0}^\infty a_n(\tilde{s}-\tilde{s}_0)^n $ which converges in a
neighborhood $\tilde{V}$ (resp., $\tilde{W}$) of $\tilde{s}_0$
in $\tilde{\omega}$ in the sense of sharp topology of $\tilde{\omega}$
(resp., in $\tilde{\Omega}$ in the sense of sharp topology of $\tilde{\Omega}$)  such
that
\begin{equation} \label{red}
f(\tilde{s})=\sum_{|n|=0}^\infty a_n(\tilde{s}-\tilde{s}_0)^n ,\;  \tilde{s}
\in \tilde{V},\; (\mbox{ resp., } \tilde{s}\in \tilde{W}).
\end{equation}
\noindent(
$(\tilde{s}-\tilde{s}_0)^n=(\tilde{s}_1-\tilde{s}_{0,1})^{n_1}...(\tilde{s}_d-
\tilde{s}_{0,d})^{n_d}$.) Actually, we will use (\ref{red}) with a representative
$(f_\varepsilon)_\varepsilon$ and $|s_\varepsilon-s_{0,\varepsilon}|=
O(\varepsilon^{-\ln r})$ (cf. the definition of the ball
$\tilde{B}(s_0,r)$.)
 It is said
that $f$ is S-real analytic in $\omega$ (resp., S-holomorphic in $\Omega$) if $f$ is S-real analytic in $s_0$
(resp.,  S-holomorphic in $s_0$) for all
$s_0\in \tilde{\omega}_c$ (resp., $s_0\in \tilde{\Omega}_c$). We denote by $\Ge_{SA}(\omega)$
(resp., $\Ge_{SH}(\omega)$) the
algebra of S-real-analytic (resp., S-holomorphic) generalized functions.
\end{definition}


$\Ge_{SA}(\omega)$, (resp., $\Ge_{SH}(\Omega)$)
 is a differential subalgebra of $\Ge(\omega)$ (resp.,
 $\Ge(\Omega)$).

Recall, a series $\sum_n a_n$ converges in an ultrametric space if
and only if its general therm $a_n$ tends to zero as $n\rightarrow \infty$.

In the complex case, the  sets of of S-holomorphic and holomorphic generalized functions  coincide.
Thus we extend the result
of \cite{ar2} where it is proved in the one dimensional case  that if   $f $ is holomorphic and sub-linear then
$f$ is S-holomorphic.
\begin{theorem}  \label{sss}
 Let $\Omega$ be an open  set in $\mathbb{C}^p$ and
$f\in\Ge(\Omega)$.
Then $f $ is holomorphic in $\Omega$ iff
$f$ is S-holomorphic in $\tilde{\Omega}$.
\end{theorem}
\bproof If $[(f_\varepsilon)_\varepsilon]$ is $S-$holomorphic, then
(as it is proved in \cite{ar2}) derivatives can be made on the
series (\ref{red}) (for $f_\varepsilon$) term by term which gives
$\bar{\partial}f=0.$ Conversely, if $\bar{\partial}f=0$ and
$\tilde{s}_0
=[(s_{0\varepsilon})_\varepsilon]\in\tilde{\Omega}_c,$
 then
by the Cauchy formula, for every $\alpha \in \mathbb{N}_0^d,$ we obtain
bounds of the form
$$|f_\varepsilon^{(\alpha)}(s_\varepsilon)|\leq C\eta^{|\alpha|}\alpha!\varepsilon^{-a},
\varepsilon<\varepsilon_0,
\mbox{  for some } a>0,
$$
for $\tilde{s}$ in a sharp neighborhood of $\tilde{s}_0.$ It follows that the corresponding Taylor
series, expressed in terms of $z$ and $\bar{z}$ converges  in the sharp topology, since
terms of $\bar{\partial}f_\varepsilon$ involving $\bar{z}$ are all equal to zero. Thus,  $f$ is $S$-analytic.
\eproof
Now we extend the notion of sub-linearity of \cite{ar2} in the $d-$dimensional case.
\begin{definition} \label{subl}
(\cite{ar2})  An $f\in \Ge(\omega)$ is sub-linear in an open set $\omega\subset
\mathbb{R}^d$ if there exists a representative $(f_\eps)_\eps$  of
$f$ with the following property: For every $\tilde{x}\in\tilde{\omega}_c$
there exist its representative $(x_\eps)_\eps$ of $\tilde{x}$, $k\in
\mathbb{R}$ and a sequence $(p_n)_n$  in $\mathbb{R}$ such that
$$\lim_{n\rightarrow \infty}(p_n+kn)=\infty \mbox{ and }
|f_\eps^{(\alpha)}(x_\eps)|=O(\eps^{p_n}),\; |\alpha|=n,  n\in\mathbb{N}_0.$$
\end{definition}
Note if the conditions of this  definition hold for some representatives of $f$ and
$\tilde{x}$, then they hold for any representatives of $f$ and
$\tilde{x}$.

The set of all sub-linear functions is a
$\mathbb{C}$-subalgebra of $\Ge(\omega)$ and
$\Ge^\infty(\omega)$ is its subalgebra.

Let $f\in\mathcal{E}'(\omega)$ and $[(f*{\phi_\eps}|_{\omega})_\eps]$ be the corresponding (embedded) generalized
function. Then it is
 a sub-linear generalized function. The same holds for $f\in\mathcal{D}'(\omega)$.

\begin{theorem}  \label{sss}

 a) Let $\Omega$ be an open  set in
$\mathbb{C}^p$ and $f\in\Ge(\Omega)$.

 Let $f$ be S-holomorphic in $\Omega$. Then $f$ is sub-linear in $\Omega\subset\mathbb{R}^{2p}$.


b) Let $\omega$ be an open  set in
$\mathbb{R}^d$,  and $f\in\Ge(\omega)$.

1. If $f$ is S-real analytic then it is sub-linear.

2. If $f $ is real analytic, then $f$ is S-real analytic.

3. If $T\in\mathcal{E}'(\omega),$ then  $[(T*\phi_\varepsilon)_\varepsilon]\in\mathcal{G}(\omega)$
 is S-real analytic.

In particular, every element of $\Ge^\infty(\omega )$
is S-real analytic but not necessarily real analytic.

\end{theorem}
\bproof
a)
 This result, in case of $\tilde{\mathbb{C}},$ is proved in \cite{ar2}.
 Let $\tilde{z}_0=[(z_{0,\eps})_\eps]$ be a  compactly supported  generalized point.
Since $|\partial^\alpha f(\tilde{z}_0)(\tilde{z}-\tilde{z}_0)^\alpha/\alpha!|
\rightarrow 0$ as $|\alpha|\rightarrow \infty,$
if $\tilde{z}$ belongs to a (sharp) neighborhood of $\tilde{z}_0,$ one  derives the sub-linearity at $\tilde{z}_0.$

b) The proof of 1. is the same as for holomorphic generalized
functions. Concerning 2., it comes from the remark that
$\mathcal{G}_A(\omega)\subset\mathcal{G}^\infty(\omega)$ and 3. So
let us prove 3.

3. Let a point
$\tilde{x}\in \tilde{\omega}_c$ be supported by $\omega''\subset\subset
\omega, $ and $\omega''\subset\subset \omega'\subset\subset
\omega$. By $T*\phi_\varepsilon^{(\alpha)}=(T*\phi_\varepsilon)^{(\alpha)}$,  we have
that there exists $k>0$ such that
 $$\sup_{t\in\omega''}
|T*\phi_\varepsilon^{(\alpha)}(t)|=O(\eps^{-k-|\alpha|}), \;\alpha \in \mathbb{N}^d.$$

With $$a_n=
[(\frac{T*\phi_\varepsilon^{(n)}(x_\eps)}{n!})_\eps],\; n\in
\mathbb{N}_0^d,$$ we have $a_{n,\eps}=O(\eps^{-k-|n|}).$ Consider the ball $\tilde{B}(\tilde{x},r)$ so that
$-\ln r =k+2$ ($|y_\varepsilon-x_\varepsilon|<\varepsilon^{k+2}$). Thus, for $\tilde{y}\in\tilde{B}(\tilde{x},r),$
$ a_n|y-x|^{|n|}$  converges to zero
in the sense of sharp topology as $|n|\rightarrow \infty$ and
$$f(\tilde{y})=\sum_n
a_n(\tilde{y}-\tilde{x})^n ,\; (|y_\eps-x_\eps|\leq \varepsilon^{k+2}),$$ converges.
Hence $f$ is an S-analytic generalized function.

 The proof of the particular case is omitted.

 \vspace{0,2cm}

In general, we can conclude that
 an S-analytic generalized function can be
irregular while  an analytic generalized function is very
regular.

{\bf Acknowledgement.} Work of the first author is partially
supported by START-project Y237 of the  Austrian Science Fund and
by the MNZZS of Serbia 144016

%


\begin {thebibliography}{13}


\bibitem{ar1} Aragona, J.,
{On existence theorems for the
    $\bar\partial$ operator on
generalized differential forms},
 Proc. London Math. Soc.  53(1986), 474--488.

\bibitem{ar2}  Aragona, J., Fernandez, R.,   Juriaans, S. O.,
{Discontinuous Colombeau
differential calculus}, Monatsh. Math. 144(2005), 13--29.


\bibitem{bia} Biagioni, H. A.,
{A Nonlinear Theory of Generalized
Functions}, Springer-Verlag, Berlin-Hedelberg-New York, 1990.

\bibitem{col} Colombeau, J. F.,
{New Generalized Functions and
Multiplications of Distributions}, North Holland, 1982.

\bibitem{co1} Colombeau, J. F.,
{Elementary Introduction in New
Generalized Functions}, North Holland, 1985.

\bibitem{coga2} Colombeau, J. F.,  Gal{\'e}, J. E.,
{Holomorphic
  generalized functions},
 J. Math. Anal. Appl.  103(1984), 117--133.



\bibitem{coga} Colombeau, J. F.,  Gal{\'e}, J. E.,
{Analytic continuation
of generalized functions}, Acta Math. Hung.\ 52(1988), 57-60.

\bibitem{gar1}  Garetto, C.,
Topological structures in Colombeau algebras:
investigation of the duals of $\mathcal{G}_c(\Omega)$, $\mathcal{G}(\Omega)$ and $\mathcal{GS}(\R^n)$,
Monatsh. Math. 146 (2005), 203-226.

\bibitem{gar2} Garetto, C., Topological structures in Colombeau algebras:
topological $\widetilde{\mathbb{C}}$-modules and duality theory,
Acta Appl. Math., 88 (2005), 81-123.

\bibitem{garhor}  Garetto, C.,  H\" ormann, G.,
Microlocal analysis of generalized functions:
pseudodifferential techniques and propagation of singularities,
Proc. Edinburgh Math. Soc. 48, (2005), 603-629.


\bibitem{gkos} Grosser, M., Kunzinger, M., Oberguggenberger, M.,
Steinbauer,R.,
{Geometric Generalized Functions with Applications to General
Relativity}, Kluwer, 2001.

\bibitem{hor}
 H\"ormander, L., The analysis of linear partial differential
operators. I.
Grundlehren der Mathematischen Wissenschaften, 256.
Springer-Verlag, Berlin, 1983.


\bibitem{horm}
 Hormander, L., {Uniqueness theorems and wave front sets for
solutions of linear differential equations with analytic
coefficients}, Comm. on pure and Appl. Math. 24 (1971) 671-704.

\bibitem{horkun}  H\" ormann, G., Kunzinger, M.,
 Microlocal properties of basic operations in Colombeau algebras,
J. Math. Anal. Appl. 261, 254-270 (2001).

\bibitem{horman}  H\" ormann, G.,  H\" older-Zygmund regularity in algebras of generalized functions,
Z. Anal. Anwendungen 23 (2004), 139-165.


\bibitem{kelsca}
{ Khelif, A.,  Scarpalezos, D.,} {Zeros of generalized holomorphic
functions}, Monach.Math, (2006), to appear.


\bibitem{mar}
{ Marti, J. A.,} $\mathcal{G}_L$-microlocal analysis of
generalized functions. Integral Transforms Spec. Funct. 17 (2006),
119--125.

\bibitem{nps} { Nedeljkov, M.,  Pilipovi\' c, S.,  Scarpalezeos, D.,},
Linear Theory of Colombeau's Generalized Functions, Addison Wesley, Longman, 1998.

\bibitem{ober 001}
{ Oberguggenberger, M.,} {Multiplication of distributions and
application to partial differential equations}. Pitman Res. Notes
Math. Ser. 259, Longman, Harlow, 1992.

\bibitem{kuob}
{ Oberguggenberger, M., Kunzinger, M.,}{Characterization of
 Colombeau generalized functions by their pointvalues}. Math. Nachr. 203 (1999), 147--157.

\bibitem{obpisca}
{ Oberguggenberger, M., Pilipovi\' c, S., Scarpalezos, D.,}
{Positivity and positive definiteness in  generalized function
algebras}. J.Math. Anal. Appl. (2006), to appear.



\bibitem{pili}
{Oberguggenberger, M., Pilipovi\' c, S., Valmorin, V.,} {Global
holomorphic representatives of holomorphic generalized functions}
Monach. Math. 2006, to appear.

\bibitem{hp}
{ Pilipovi\' c, S.,} {Generalized hyperfunctions and algebra of
megafunctions}, Tokyo J. Math.  28 (2005),  1--12.

\bibitem{dim} { Scarpalezos, D.,}
{Topologies dans les espaces de nouvelles fonctions g{\'e}n{\'e}ralis{\'e}es de
Colombeau. {$\widetilde{\mathbb C}$}--modules
     topologiques},
 {Preprint ser. Universit{\'e} Paris 7},
     {Universit{\'e} Paris 7} {1992}.

\bibitem{dimbec } { Scarpalezos, D.,}
{Some remarks on functoriality of Colombeau's construction; topological
and microlocal aspects and applications.} Generalized functions---linear and nonlinear
problems (Novi Sad, 1996, Integral Transform. Spec. Funct. 6 (1998), 295--307.

\bibitem{dimint } { Scarpalezos, D.,}
{Colombeau's generalized functions: topological structures; microlocal properties.
 A simplified point of view, I}, Bull. Cl. Sci. Math. Nat. Sci. Math., (2000), 89--114.

\bibitem{sch} { Schwartz, L.,} Theorie des distributiones, Hermann, Paris, 1950.

\bibitem{valm}
{Valmorin, V.,}
{ Vanishing theorems in Colombeau  generalized functions}, 
Canad. Math. Bull., to appear.
\end{thebibliography}
\end{document}